\documentclass[12pt]{article}
\usepackage{color,latexsym,amsmath,amssymb,path}
\usepackage{enumitem}

\newtheorem{theorem}{Theorem}[section]

\newtheorem{corollary}[theorem]{Corollary}

\newcommand{\Proof}[1]
        {
        \noindent
        \emph{Proof #1.}~
        }
\newsavebox{\smallProofsym}                     
\savebox{\smallProofsym}                        %
        {
        \begin{picture}(7,7)                    %
        \put(0,0){\framebox(6,6){}}             %
        \put(0,2){\framebox(4,4){}}             %
        \end{picture}                           %
        }                                       %
\newcommand{\smalleop}[1]
        {
        \mbox{} \hfill #1~~\usebox{\smallProofsym}\!\!\!\!\!\!\
        }

\usepackage{graphicx,psfrag}

\newcommand{\placefig}[2]
        {\includegraphics[width=#2]{#1.eps}}
\usepackage{dsfont}

\newcommand{\RR}{\ensuremath{\mathbb R}}

\newcommand{\pts}{\mathcal P}
\newcommand{\ptss}{\mathcal V}
\newcommand{\curves}{\mathcal C}

\def\eps{{\varepsilon}}

\begin{document}
\pagenumbering{arabic}

\title{Distinct distances on two lines\thanks{%
Work on this paper by the first two authors was partially
supported by Grant 338/09 from the Israel Science Fund
and by the Israeli Centers of Research Excellence (I-CORE)
program (Center No.~4/11). Work by Micha Sharir was also
supported by NSF Grant CCF-08-30272 and
by the Hermann Minkowski-MINERVA Center for Geometry
at Tel Aviv University. Work by J\'ozsef Solymosi was supported by NSERC, ERC-AdG 321104, and OTKA NK 104183 grants. }}


\author{
Micha Sharir\thanks{%
School of Computer Science, Tel Aviv University,
Tel Aviv 69978, Israel.
{\sl michas@tau.ac.il} }
\and
Adam Sheffer\thanks{%
School of Computer Science, Tel Aviv University,
Tel Aviv 69978, Israel. 
{\sl sheffera@tau.ac.il}}
\and
J\'ozsef Solymosi\thanks{%
Department of Mathematics,
University of British Columbia,
Vancouver, BC, V6T 1Z4, Canada. 
{\sl solymosi@math.ubc.ca}} }

\maketitle


\begin{abstract}
Let $\pts_1$ and $\pts_2$ be two finite sets of points in the plane, so
that $\pts_1$ is contained in a line $\ell_1$, $\pts_2$ is contained in
a line $\ell_2$, and $\ell_1$ and $\ell_2$ are neither parallel nor
orthogonal. Then the number of distinct distances determined by the
pairs of $\pts_1 \times \pts_2$ is
\[ \Omega\left( \min\left\{|\pts_1|^{2/3}|\pts_2|^{2/3},
|\pts_1|^2,|\pts_2|^2 \right\}\right). \]
In particular, if $|\pts_1|=|\pts_2|=m$, then the number of these distinct
distances is $\Omega(m^{4/3})$, improving upon the previous bound
$\Omega(m^{5/4})$ of Elekes \cite{Elek99}.
\end{abstract}


\section{Introduction}

Given a set $\pts$ of $m$ points in $\RR^2$, let $D(\pts)$ denote the number
of distinct distances that are determined by pairs of points from $\pts$.
Let $D(m) = \min_{|\pts|=m}D(\pts)$; that is, $D(m)$ is the minimum number of
distinct distances that any set of $m$ points in $\RR^2$ must always determine.
In his celebrated 1946 paper \cite{erd46}, Erd\H os derived the bound
$D(m) = O(m/\sqrt{\log m})$.
For the celebrations of his 80'th birthday, Erd\H os compiled a survey of
his favorite contributions to mathematics \cite{erd96}, in which he wrote
``My most striking contribution to geometry is, no doubt, my problem on the
number of distinct distances. This can be found in many of my papers on
combinatorial and geometric problems".  Recently, after 65 years and a series
of increasingly larger lower bounds (comprehensively described in the book
\cite{GIS11}), Guth and Katz \cite{GK11} provided an almost matching lower
bound $D(m) = \Omega(m/\log m)$. For this, Guth and Katz developed several
novel techniques, relying on tools from algebraic geometry.

While the problem of obtaining the asymptotic value of $D(m)$ is almost settled,
many other variants of the distinct distances problem are still widely open.
For example, see \cite{Dum06,NPPZ11} regarding the conjecture that any $m$
points in convex position in the plane determine at least $\lfloor m/2\rfloor$
distinct distances, and \cite{SV08} for a study of the minimum number of
distinct distances in higher dimensions.

\begin{figure}[h]
\centerline{\placefig{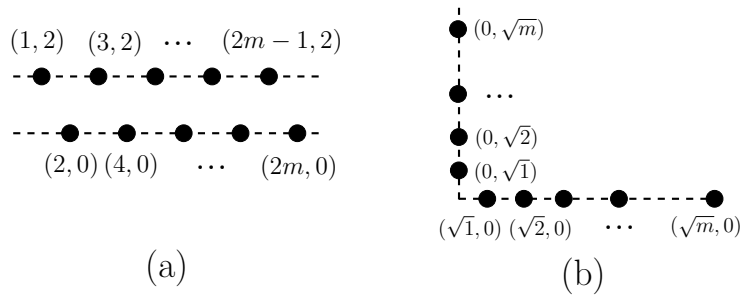}{0.64\textwidth}}
\vspace{-1mm}

\caption{\small \sf (a) Two parallel lines with $D(\pts_1,\pts_2) = \Theta(m)$. 
(b) Two orthogonal lines with $D(\pts_1,\pts_2) = \Theta(m)$.}
\label{fi:quadruple}
\vspace{-2mm}
\end{figure}

In this paper we consider the following variant of the distinct distances
problem in the plane.  Let $\pts_1$ and $\pts_2$ be two finite sets of points,
such that all the points of $\pts_1$ lie on a line $\ell_1$, and all the points
of $\pts_2$ lie on a line $\ell_2$. 
Let $D(\pts_1,\pts_2)$ denote the number of distinct
distances between the points of $\pts_1$ and $\pts_2$, i.e.,
$$D(\pts_1,\pts_2)=\biggl|\{\text{dist}(p,q) \;| \; p\in \pts_1, q\in \pts_2\}\biggr|.$$  
Consider first the ``balanced'' case, where $|\pts_1|=|\pts_2|=m$. When the two 
lines are parallel or orthogonal, the points can be arranged such that
$D(\pts_1,\pts_2) = \Theta(m)$; for example, see Figure \ref{fi:quadruple}.
Purdy conjectured that if the lines are neither parallel nor orthogonal
then $D(\pts_1,\pts_2) = \omega(m)$  (e.g., see \cite[Section 5.5]{BMP05}).
Elekes and R\'onyai \cite{ER00} proved that the number of distinct distances
in such a scenario is indeed superlinear. They did not give an explicit bound,
but a brief analysis of their proof shows that $D(\pts_1,\pts_2) = \Omega(m^{1+\delta})$
for some $\delta>0$. Elekes \cite{Elek99} derived the improved bound
$D(\pts_1,\pts_2) = \Omega(m^{5/4})$ (when the lines are neither parallel nor
orthogonal) and gave a construction, reminiscent of the one by Erd{\H o}s~\cite{erd46},
with $D(\pts_1,\pts_2) = O\left(m^2/\sqrt{\log m}\right)$, in which the angle
between the two lines is $\pi/3$. Previously, these were the best known bounds
for $D(\pts_1,\pts_2)$ for the balanced case. The unbalanced case, where 
$|\pts_1|\ne |\pts_2|$, has recently been studied by Schwartz, Solymosi, 
and de Zeeuw \cite{SSZ12}, who have shown, among several other related results, 
that the number of distinct distances remains superlinear when 
$|\pts_1| = m^{1/2 + \varepsilon}$ and $|\pts_2| = m$, for any $\varepsilon>0$.

In this paper we derive the following result, for point sets $\pts_1,\pts_2$
of arbitrary (possibly different) cardinalities.

\begin{theorem} \label{th:main}
Let $\pts_1$ and $\pts_2$ be two sets of points in $\RR^2$ of cardinalities $n$
and $m$, respectively, such that the points of $\pts_1$ all lie on a line $\ell_1$,
the points of $\pts_2$ all lie on a line $\ell_2$, and the two lines are neither
parallel nor orthogonal. Then the number of distinct distances between $\pts_1$
and $\pts_2$ is
\[D(\pts_1,\pts_2)=\Omega\left( \min\left\{n^{2/3}m^{2/3},n^2,m^2 \right\}\right). \]
\end{theorem}

Theorem \ref{th:main} immediately implies the following improved lower bound for
the balanced case.

\begin{corollary} \label{co:bal}
Let $\pts_1$ and $\pts_2$ be two sets of points in $\RR^2$, each of cardinality
$m$, such that the points of $\pts_1$ all lie on a line $\ell_1$, the points of
$\pts_2$ all lie on a line $\ell_2$, and the two lines are neither parallel nor
orthogonal. Then $D(\pts_1,\pts_2)=\Omega\left( m^{4/3}\right).$
\end{corollary}

Note that Theorem~\ref{th:main} also implies the result of \cite{SSZ12},
mentioned earlier, and slightly strengthens it, by providing the explicit 
lower bound $\Omega(m^{1+2\eps/3})$.

Even with the improved lower bound in Corollary~\ref{co:bal} (over the lower
bound in \cite{Elek99}), there is still a considerable gap to the near-quadratic
upper bound in \cite{Elek99}, and the prevailing belief is that the correct lower
bound is indeed close to quadratic.

To obtain the improved bound, we use a double counting argument, applied to
the number of quadruples $(a,p,b,q)$ of points, with $a,b\in \pts_1$ and
$p,q\in\pts_2$, that satisfy $|ap|=|bq|$. The argument is similar to
the one in the reduction devised by Elekes and presented in Elekes and
Sharir \cite{ES11}. For this double counting we use the same lower bound
analysis as in \cite{ES11}, but replace the upper bound analysis by a
considerably simpler one, in which the problem is reduced to that of bounding
the number of incidences between certain points and hyperbolas in the plane.
(In contrast, the original reduction in \cite{ES11} is to incidences between
points and lines in $\RR^3$.)



\section{The proof of Theorem \ref{th:main}} \label{sec:proof}

Without loss of generality, we may assume that the points of $\pts_1$ are on one
side of the intersection point $\ell_1\cap\ell_2$.  Otherwise, we can partition $\pts_1$
into two subsets by splitting $\ell_1$ at $\ell_1\cap\ell_2$, and remove the
subset that yields fewer distinct distances with the points of $\pts_2$.
At worst, this halves the number of distinct distances between the pairs in
$\pts_1\times\pts_2$. For the same reason, we may assume that the points of
$\pts_2$ are also on one side of $\ell_1\cap\ell_2$. Furthermore, without
loss of generality, we may assume that $n = |\pts_1| \ge m = |\pts_2|$.

We rotate, translate, and possibly reflect the original plane, so that the
origin $o$ is $\ell_1 \cap \ell_2$, $\ell_1$ is the $x$-axis, the points of
$\pts_1$ lie on the positive side of $o$, and the points of $\pts_2$ lie
above $\ell_1$. We denote the angle between the two lines by $\alpha$. 
Since the two lines $\ell_1$ and $\ell_2$ are neither parallel nor orthogonal,
we have $\alpha\neq 0,\pi/2$.
We will also assume that $o\notin \pts_1\cup\pts_2$, because the presence of $o$
in either set can generate at most $O(m+n)$ distinct distances.

\begin{figure}[h]
\centerline{\placefig{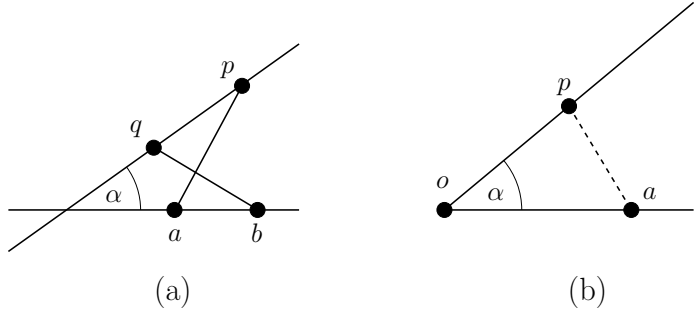}{0.6\textwidth}}
\vspace{-1mm}

\caption{\small \sf (a) A quadruple $(a,p,b,q)$ in $Q$.
(b) By the law of cosines, we have $|ap|^2 = |oa|^2 + |op|^2 -2|oa| |op|\cos \alpha$.}
\label{fi:cosine}
\vspace{-2mm}
\end{figure}

We begin with a variant of the first part of the reduction from \cite{ES11}.
We set $x=D(\pts_1,\pts_2)$ and denote the $x$ distinct distances in
$\pts_1\times\pts_2$ as $\delta_1,\ldots, \delta_x$. For a pair of points $u$
and $v$, we denote by $|uv|$ the (Euclidean) distance between $u$ and $v$.
Let $Q$ be the set of quadruples $(a,p,b,q)$, where $a,b \in \pts_1$ and
$p,q \in \pts_2$, such that $|ap|=|bq|>0$ and $ap \neq bq$ (the two segments
are allowed to share at most one endpoint); see Figure \ref{fi:cosine}(a).
The quadruples are ordered, so that $(a,p,b,q)$ and $(b,q,a,p)$ are considered
as two distinct elements of $Q$.

Let $E_i = \{(a,p)\in \pts_1 \times \pts_2 \mid |ap|=\delta_i \}$, for
$i=1,\ldots,x$. Using the Cauchy-Schwarz inequality we have,
\begin{equation} \label{eq:lowerQbi}
|Q| = 2\sum_{i=1}^x\binom{|E_i|}{2} \ge
\sum_{i=1}^x \left(|E_i|-1 \right)^2 \ge \frac{1}{x}\left(\sum_{i=1}^x (|E_i|-1) \right)^2 =
\frac{\left(mn -x\right)^2}{x}.
\end{equation}

In the remainder of the proof we derive an upper bound on $|Q|$, showing that
$|Q|=O(m^{4/3}n^{4/3} + n^2)$. Combined with \eqref{eq:lowerQbi} this yields
the lower bound asserted in the theorem (under the assumption $n\ge m$).

To obtain this upper bound, we re-interpret $|Q|$ as an incidence count between
certain points and hyperbolas in a suitable parametric plane, and then use standard machinery to
bound this count. This replaces (and greatly simplifies) the second part of Elekes's
reduction, where $|Q|$ is interpreted as an incidence count between points and lines
in $\RR^3$.

Let us consider a quadruple $(a,p,b,q)\in (\pts_1 \times \pts_2)^2$.
By the law of cosines, we have $|ap|^2 = |oa|^2 + |op|^2 -2|oa| |op|\cos \alpha$
and $|bq|^2 = |ob|^2 + |oq|^2 -2|ob| |oq|\cos \alpha$ (see Figure \ref{fi:cosine}(b)).
Thus, the quadruple $(a,p,b,q)$ is in $Q$ if and only if
\[ |oa|^2 + |op|^2 -2|oa| |op|\cos \alpha = |ob|^2 + |oq|^2 -2|ob| |oq|\cos \alpha. \]
We 
represent each point $u$ of $\pts_1$ or of $\pts_2$ by its 
distance $|ou|$ from the origin $o.$ Each of $\pts_1$, $\pts_2$
is contained in a ray (with initial point $o$) of the respective line $\ell_1$, $\ell_2$, we may assume
that all these distances are all distinct in $\pts_1$ and are all distinct in $\pts_2$. 

In what follows, we will use $u$ to denote both the point and its distance $|ou|$ from $o$.
Using the notation $s=\cos\alpha$, the above condition can be written as
\begin{equation} \label{eq:condition}
a^2-b^2 + p^2-q^2 -2s(ap-bq) = 0 ,
\end{equation}
where $s\ne 0,1$.

Let $\ptss_1$ and $\ptss_2$ denote the sets of ordered distinct pairs of $\pts_1$ and of
$\pts_2$, respectively. That is,
\[\ptss_i = \pts_i\times\pts_i\setminus \{(x,x) \mid x\in\pts_i\}\;\quad \text{for} \; i=1,2.\]

For every  $i\in\{1,2\}$ and $(a,b) \in \ptss_i$,
there is a curve $\gamma_{a,b}^{(i)}$, corresponding to the pair $(a,b)$, given by the equation
\begin{equation}\label{eq:hyperDef}
a^2-b^2 + x^2-y^2 -2s(ax - by) = 0.
\end{equation}
This can also be written as
$$
(x-sa)^2-(y-sb)^2 = (1-s^2)(b^2-a^2) .
$$
Since $s\neq 1$ and $a\ne b$, the curve $\gamma_{a,b}^{(i)}$ is a hyperbola.
Moreover, since $s\ne 0$,
all the hyperbolas $\gamma_{a,b}^{(1)}$ are distinct, and so are
all the hyperbolas $\gamma_{a,b}^{(2)}$.
Let $\curves_1$ denote the set of the $m(m-1)$
hyperbolas, $\curves_1=\{\gamma_{a,b}^{(1)}\; | (a,b) \in \ptss_1\}$.
By construction, the hyperbola $\gamma_{a,b}^{(1)} \in \curves_1$ is incident to the
point $(p,q) \in \ptss_2$ if (and only if) $a,b,p,q$ satisfy the condition in
(\ref{eq:condition}). That is, the number of quadruples $(a,p,b,q)$ in $Q$
for which $a\ne b$ and $p\ne q$ is at most the number of
point-hyperbola incidences between $\ptss_2$ and $\curves_1$. 
The number of missing elements of $Q$, where either $a=b$ or $p=q$, 
has the trivial upper bound $4mn$, as is easily verified.

Every pair of hyperbolas from $\curves_1$ intersect in at most two points 
in the real plane. Therefore, for every pair of hyperbolas from $\curves_1$ there
are at most two points of $\ptss_2$ that are incident to both hyperbolas.
The roles of $\ptss_1$ and $\ptss_2$ in (\ref{eq:hyperDef}) are symmetric. 
In particular, every pair of hyperbolas from $\curves_2=\{\gamma_{a,b}^{(2)} \mid (a,b) \in \ptss_2\}$ 
intersect in at most two points too. This implies that there are at most two 
hyperbolas of $\curves_1$ that pass through a given pair of points in $\ptss_2$.

We can therefore use the following result of Pach and Sharir \cite{PS98}.

\begin{theorem} \label{th:PS} {\bf (Pach and Sharir \cite{PS98})}
Consider a point set $\pts$, a set of curves $\cal C$, and a constant
positive integer $s$, such that 
\begin{description}
\item[(i)] for every pair of points of $\pts$ there
are at most $s$ curves of $\cal C$ that are incident to both points, and
\item[(ii)] every pair of curves of $\cal C$ intersect in at most $s$
points. 
\end{description}
Then the number of incidences between $\pts$ and $\cal C$ is
at most
$\displaystyle c(s)(|\pts|^{2/3}|{\cal C}|^{2/3}+|\pts|+|{\cal C}|)$,
where $c(s)$ is a constant that depends on $s$.
\end{theorem}

We apply Theorem \ref{th:PS} to $\ptss_2$ and $\curves_1$, with
$s=2$, and obtain (adding the bound on the missing pairs in $Q$)
$$
|Q|=O(m^{4/3}n^{4/3}+m^2+n^2+mn)
=O(m^{4/3}n^{4/3}+n^2)
$$
(recalling that we assume $n\ge m$), as desired. As already noted,
combining this bound with \eqref{eq:lowerQbi} implies
\[ \frac{\left(mn -x\right)^2}{x} = O(m^{4/3}n^{4/3} + n^2), \]
or, as is easily checked,
\[ x = \Omega\left(\min\left\{m^{2/3}n^{2/3},m^2\right\} \right).\]
Combining this bound with its symmetric version when $m\ge n$ yields
the bound asserted in the theorem. \smalleop{}
\vspace{2mm}

%

\noindent {\bf Acknowledgements.}
The second author would like to thank Andrew Suk for introducing this problem to him.


\end{document}